\documentstyle[proceedings]{article} 
\title{Extremal Optimization: Methods derived from Co-Evolution} 
\author{{\bf Stefan Boettcher\thanks{~~e-mail: sboettc@emory.edu}} \\  
Physics Department\\  
Emory University\\ 
Atlanta, GA 30322 \\ 
\And 
{\bf Allon G. Percus\thanks{~~e-mail: percus@lanl.gov}}  \\ 
Computer Research and Applications Group (CIC-3)\\
Los Alamos National Laboratory \\ 
Los Alamos, NM 87545 \\ 
} 
\begin{document} 
\maketitle 
\begin{abstract} 
We describe a general-purpose method for finding high-quality solutions
to hard optimization problems, inspired by self-organized
critical models of co-evolution such as the Bak-Sneppen model.
The method, called {\em Extremal Optimization\/}, successively eliminates
extremely undesirable components of sub-optimal solutions, rather than 
``breeding'' better components. In contrast to {\em Genetic Algorithms\/} which
operate on an entire ``gene-pool'' of possible solutions, Extremal
Optimization improves on a single candidate solution by treating each
of its components as species co-evolving according to Darwinian
principles. Unlike {\em Simulated Annealing,\/} its non-equilibrium
approach effects an algorithm requiring few parameters to tune.
With only one adjustable parameter, its performance  proves competitive
with,  and often superior to, more elaborate stochastic optimization
procedures.  We demonstrate it here on two classic hard optimization
problems: graph  partitioning and the traveling salesman problem.
\end{abstract} 
 
\section{Natural Emergence of Optimized Configurations}

Every day, enormous efforts are devoted to organizing the supply and
demand of limited resources, so as to optimize
their utility.  Examples include the supply of foods and services
to consumers, the scheduling of a transportation fleet, or the flow of
information in communication networks within society or within a
parallel computer.  By contrast, {\em without\/} any intelligent
organizing facility, many natural systems have evolved into amazingly
complex structures that optimize the utilization of resources in
surprisingly sophisticated ways (Bak 1996). For instance, driven
merely by sunlight, biological evolution has developed efficient and
strongly interdependent networks in which resources rarely go to
waste.  Even the inanimate morphology of natural landscapes exhibits
patterns far from random that often seem to serve a purpose, such as
the efficient drainage of water (Rodriguez-Iturbe 1997). The physical properties
of these fractal patterns have aroused the interest of statistical
physicists in recent times (Mandelbrot 1983).

Natural systems that exhibit such self-organizing qualities often
possess common features: they generally consist of a large number of
strongly coupled entities with very similar properties (like
species in biological evolution, despite their apparent differences).
Hence, they permit a statistical description at some coarse level. An
external resource (such as sunlight) drives the system which then
takes its direction purely by chance. If we were to rerun evolution,
there may not be trees and elephants, say, but other complex structures. Like
flowing water breaking through the weakest of all barriers in its
wake, species are coupled in a global comparative process that
persistently washes away the least fit. In this process, unlikely but
highly adapted structures surface inadvertently, as Darwin observed
(Darwin 1859). Optimal adaptation thus emerges naturally, without
divine intervention, from the dynamics through a selection {\em
against\/} the extremely ``bad''.  In fact, this process prevents the
inflexibility that would inevitably arise in a controlled breeding of the
``good''.

Certain models relying on extremal processes have been proposed 
to explain self-organizing systems in nature (Paczuski 1996). In
particular, the Bak-Sneppen model of biological evolution is based on
this principle (Bak 1993, Sneppen 1995). It is happily devoid of any specificity
about the nature of interactions between species, yet produces 
salient nontrivial features of paleontological data such as
broadly distributed lifetimes of species, large extinction events, and
punctuated equilibrium (Gould 1977). 

\begin{figure}
\vskip 4.40truein
\includegraphics{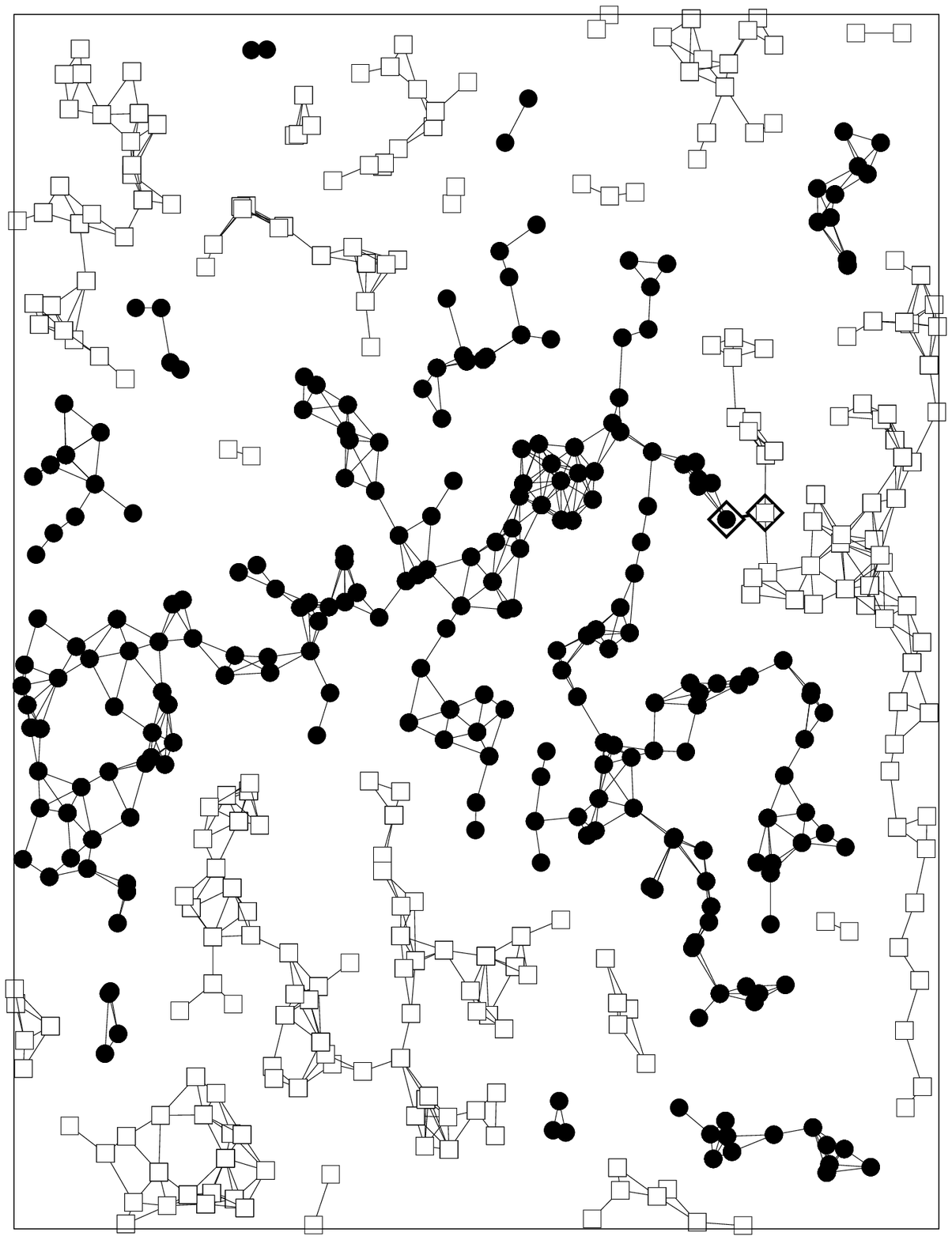}
\smallskip
\includegraphics{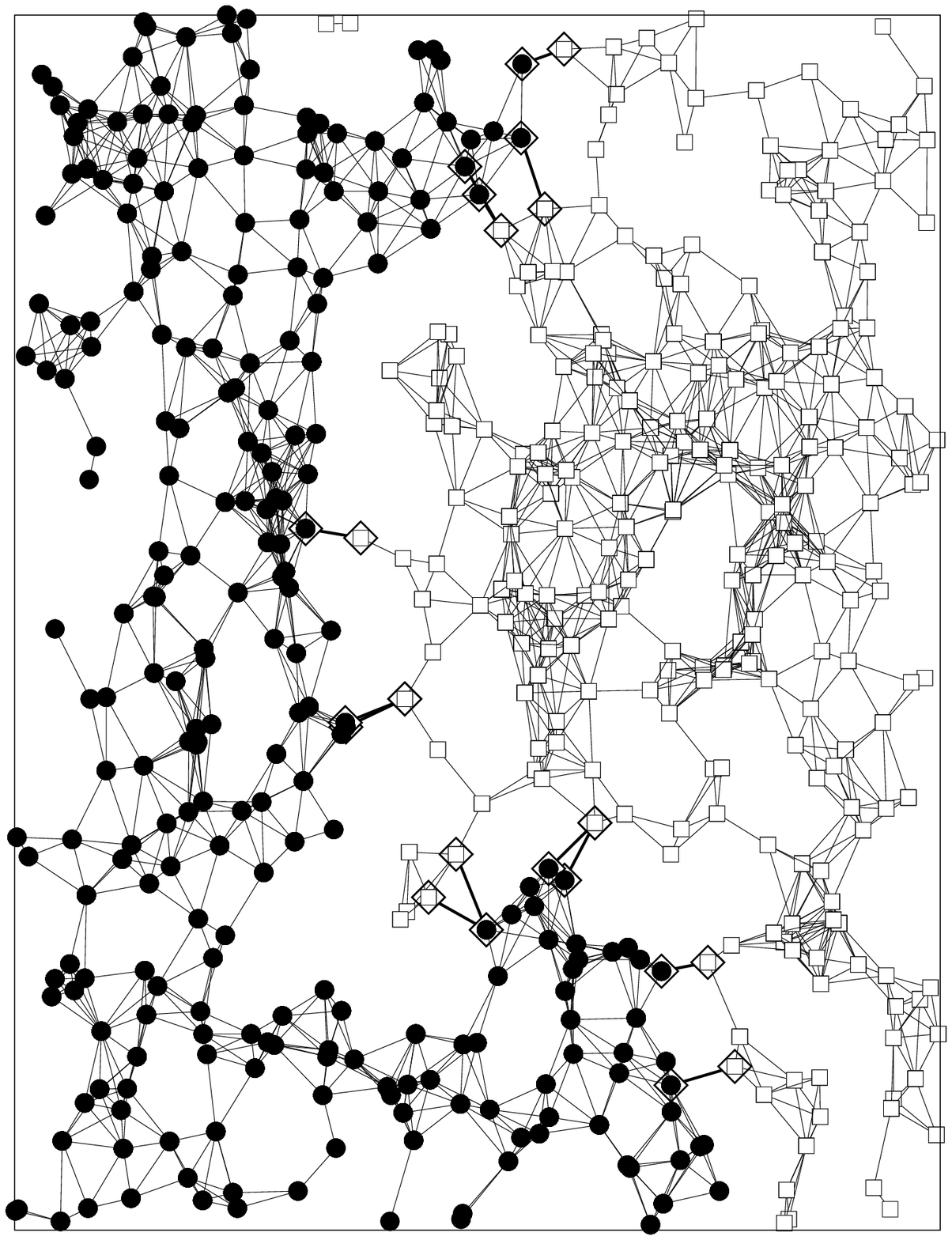}
\caption{Two random geometric graphs, $N=500$, with connectivities
$\alpha\approx4$ (top) and $\alpha\approx8$ (bottom) in an optimized
configuration found by EO. At $\alpha=4$ the
graph barely percolates, with only one ``bad'' edge connecting the set of
250 round points with the set of 250 square points (diamonds show the two
ends of the edge), thus $m_{\rm opt}=1$. For the
denser graph on the bottom, EO obtained the cutsize $m_{\rm opt}=13$.
}
\label{geograph}
\end{figure}

Species in the Bak-Sneppen model
are located on the sites of a lattice, and each is represented by a value
between 0 and 1 indicating its ``fitness''.  At each update step, the
smallest value (representing the worst adapted species) is discarded
and replaced with a new value drawn randomly from a flat distribution on
$[0,1]$.  Without any interactions, all the fitnesses in the system
would eventually become 1.  But obvious interdependencies between species
provide constraints for balancing the system's overall fitness with that of
its members: the change in fitness of one species impacts the fitness
of an interrelated species.  Therefore, at each update step in the
Bak-Sneppen model, the fitness values on the sites {\em neighboring\/}
the smallest value are replaced with new random numbers as well.  No
explicit definition is given of the mechanism by which these neighboring
species are related.  Yet after
a certain number of updates, the system organizes itself into a
highly correlated state known as self-organized criticality (SOC)
(Bak 1987). In that state, almost all species have reached a fitness
above a certain threshold.  But these also species possess what is called
punctuated equilibrium (Gould 1977): since one's weakened neighbor can undermine
one's own fitness, co-evolutionary activity gives rise to chain reactions.
Fluctuations that rearrange the fitness of many species occur routinely.
These fluctuations can be of the scale of the system itself, making any
possible configuration accessible.

In the Bak-Sneppen model, the high degree of adaptation of most
species is obtained by the elimination of badly adapted ones instead
of a particular ``engineering'' of better ones.  While such dynamics
might not lead to as optimal a solution as could be engineered in 
specific circumstances, it provides near-optimal solutions with a
high degree of latency for a rapid adaptation response to changes
in the resources that drive the system.

In the following we will describe an optimization method inspired 
by these insights (Boettcher, submitted, and Boettcher, to appear), called
{\em extremal optimization\/}, and study its performance for graph
partitioning and the traveling salesman problem.

\section{Extremal Optimization and Graph Partitioning}

In graph (bi-)partitioning, we are given a set of $N$ points, where
$N$ is even, and ``edges'' connecting certain pairs of
points. The problem is to partition the points into two
equal subsets, each of size $N/2$, with a minimal number of edges
cutting across the partition.  (Call the number of these edges
the ``cutsize'' $m$, and the optimal cutsize $m=m_{\rm opt}$.)
The points themselves could, for instance, be associated with
positions in the unit square. A ``geometric'' graph of
average connectivity $\alpha$ would then be formed by connecting any two
points within Euclidean distance $d$, where $N\pi d^2=\alpha$ (see
Fig.~\ref{geograph}).  Constraining the partitioned subsets to be of
fixed (equal) size makes the solution to the problem particularly
difficult.  This geometric problem resembles those found in
VLSI design, concerning the optimal partitioning of gates between 
integrated circuits (Dunlop 1985).

Graph partitioning is an {\em NP-hard\/} optimization
problem (Garey 1979): it is believed that for large $N$ the
number of steps necessary for an algorithm to find the {\em exact\/}
optimum must, in general, grow faster than any polynomial in $N$. In
practice, however, the goal is usually to find near-optimal
solutions quickly. Special-purpose heuristics to find approximate
solutions to specific NP-hard problems abound (Alpert 1995, Johnson 1997).
Alternatively, general-purpose
optimization approaches based on stochastic procedures have been
proposed, most notably {\em simulated annealing\/} (Kirkpatrick 1983, 
{\v C}erny 1985) and {\em genetic algorithms\/} (Holland 1975). These methods, 
although slower, are applicable to
problems for which no specialized heuristic exists.  Extremal
optimization (EO) falls into
the latter category, adaptable to a wide range of combinatorial
optimizations problems rather than crafted for a specific application.

In close analogy to the Bak-Sneppen model of SOC, the EO algorithm
proceeds as follows for the case of graph bi-partitioning:
\begin{enumerate}
\item Initially, partition the $N$ points at will into two equal subsets.
\item Rank each point $i$ according to its fitness,
$\lambda_i=g_i/(g_i+b_i)$, where $g_i$ is the number of (good) edges
connecting $i$ to points within the same subset, and $b_i$ is the number 
of (bad) edges connecting $i$ to the other subset. If point $i$ has no 
connections at all ($g_i=b_i=0$), let $\lambda_i=1$.  
\item Pick the least fit point, {\em i.e.\/}, the point (from either subset)
with the smallest $\lambda_i\in [0,1]$.  Pick a second point at random
from the other subset, and interchange these two points so that each
one is in the opposite subset from where it started.
\item Repeat at (2) for a preset number of times [assume $O(N)$ updates].
\end{enumerate}
\noindent The result of an EO run is defined as the best (minimum cutsize)
configuration seen so far.  All that is necessary to keep track of, then,
is the current configuration and the best so far.  

\begin{figure}
\vskip 2.2truein
\includegraphics{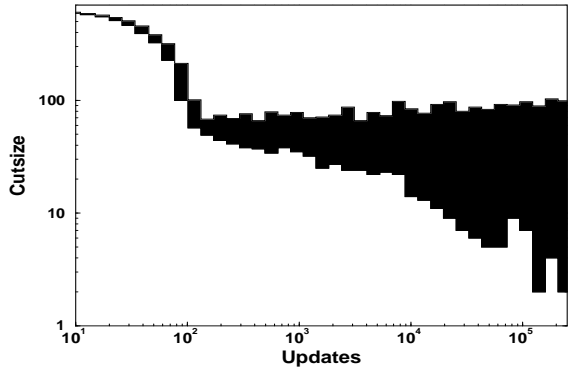}
\caption{Evolution of the cutsize during an extremal optimization run on an
$N=500$ geometric graph with $\alpha=5$ (see Fig.~\protect\ref{geograph}). The
shaded area 
marks the range of cutsizes explored in the respective time bins. The 
best cutsize ever found is 2, which is visited repeatedly in this run.
In contrast to simulated annealing, which has large fluctuations in 
early stages of the run and then converges much later, extremal optimization
quickly approaches a stage where broadly distributed fluctuations allow
it to probe many local optima. In this run, a random initial partition was
used, and the runtime on a 200MHz Pentium was 9sec.
}
\label{runtime} 
\end{figure}

EO, like simulated annealing (SA) and genetic algorithms (GA), is
inspired by  observations of physical systems [for a comparison of SA and
GA, see e. g. (de Groot 1991)]. However, SA emulates
the behavior of frustrated systems in thermal equilibrium: if one
couples such a system to a heat bath of adjustable temperature, by
cooling the system slowly one may come close to attaining a state of
minimal energy.  SA accepts or rejects local changes to a
configuration according to the  Metropolis algorithm (Metropolis 1953) at a
given  temperature, enforcing equilibrium dynamics (``detailed
balance'') and requiring a carefully tuned ``temperature schedule''.
In contrast, EO takes the system far from equilibrium: it applies no
decision criteria, and all new configurations are accepted
indiscriminately.  It may appear that EO's results would resemble an
ineffective random search. But in fact, by persistent selection
against the worst fitnesses, one quickly approaches near-optimal
solutions. At the same time, significant  fluctuations still remain at
late run-times (unlike in SA), crossing sizable barriers to access new
regions in configuration space, as shown in Fig.~\ref{runtime}. EO and
genetic  algorithms are equally contrasted.  GAs keep track of
entire ``gene pools'' of solutions from which to select and
``breed'' an improved generation of global approximations.   By
comparison, EO operates only with local updates on a single copy of
the system, with improvements achieved instead by elimination of the bad.

Further improvements may be obtained through a slight modification of the EO
procedure.  Step (2) of the algorithm establishes a fitness rank for all
points, going from rank $n=1$ for the worst fitness $\lambda$ to rank $n=N$
for the best.  (For points with degenerate values of $\lambda$, the
ranks may be assigned in random order.)  Now relax step (3) so that the
points to be interchanged are both chosen from a probability distribution
over the rank order: from each subset, we pick a point having rank $n$
with probability $P(n)\propto n^{-\tau},~1\leq n\leq N$.  The choice of a
power-law distribution for $P(n)$ ensures that no regime of fitness gets
excluded from further evolution, since $P(n)$ varies in a gradual,
scale-free manner over rank.  Universally, for a wide range of graphs, we
obtain best results for $\tau\approx 1.2-1.6$.  What is the physical
meaning of an optimal value for $\tau$?  If $\tau$ is too small, we often
dislodge already well-adapted points of high rank: ``good'' results get
destroyed too frequently and the progress of the search becomes
undirected.  On the other hand, if $\tau$ is too large, the process
approaches a deterministic local search and gets stuck near a local
optimum of poor quality.  At the optimal value of $\tau$, the more fit
components of the solution are allowed to survive, without the search
being too narrow.  Our numerical studies have indicated that the best
choice for $\tau$ is  closely related to a transition from ergodic to
non-ergodic behavior, with optimal performance of EO obtained near
the edge of ergodicity.

To evaluate EO, we tested the algorithm on a testbed of well-studied large
graphs\footnote{These instances are available at\newline
http://userwww.service.emory.edu/\~{}sboettc/graphs.html} discussed in 
(Hendrickson 1996, Merz 1998).
Table~\ref{tab1} summarizes EO's results on these, using 30 runs of at most
$200N$ update steps (in several cases far fewer were necessary; see below).
On the first four large graphs, SA's performance is extremely poor; we
therefore substitute results given in (Hendrickson 1996) using a variety of
specialized heuristics.  EO significantly improves upon these cutsizes,
though at longer runtimes.  The best results to date on the graphs are
due to various GAs (Merz 1998).  EO reproduces all of these cutsizes,
displaying an increasing runtime advantage as $N$ increases.  On the
final four graphs, for which no GA results were available, EO matches or
dramatically improves upon SA's cutsizes.  And although increasing $\alpha$ 
generally slows down
EO and speeds up SA, EO's runtime is still nearly competitive with SA's
on the high-connectivity {\em Nasa\/} graphs.

Several factors account for EO's speed.  First of all, in step (1) we employ
a simple ``greedy'' start to form the initial partition, clustering
connected points into the same partition from a random seed.  This helps
EO to succeed rapidly.  By contrast, greedy initialization improves
the performance of SA only for the smallest
and sparsest graphs.  Second of all, in step (2) we use a stochastic sorting
process to accelerate the algorithm.  At each update step, instead of
perfectly ordering the fitnesses $\lambda_i$, we arrange them on an ordered
binary tree called a ``heap''.  We then select members from the heap such
that {\em on average\/}, the actual rank selection approximates
$P(n)\sim n^{-\tau}$.  This stochastic rank sorting introduces a runtime
factor of only $\alpha\log{N}$ per update step.  Finally, EO requires significantly
fewer update steps (Fig.~\ref{runtime}) than, say, a complete SA temperature
schedule.  The quality of our large $N$ results confirms that $O(N)$ update
steps are indeed sufficient for convergence.  In the case of the {\em Nasa\/}
graphs, only $30N$ update steps (rather than the full $200N$) were in fact
required for EO to reach its best results, and in the case of the
{\em Brack2\/} graph, only $2N$ steps were required.

\begin{table}[t]
\caption{Best cutsizes and runtimes for our testbed of graphs.
EO and SA results are from our runs (SA parameters
as determined by Johnson {\em et al.\/} (Johnson 1989)), using a 200MHz
Pentium.  GA results are from Merz and Freisleben (Merz 1998), using a 300MHz
Pentium.  Comparison data for three of the large graphs are due to results
from heuristics by Hendrickson (Hendrickson 1996), using a 50MHz Sparc20.
}
\begin{center}
\begin{tabular}{lr@{\ }lr@{\ }lr@{\ }l}       %|lr@{\ }lr@{\ }lr@{\ }l}
Graph & \multicolumn{2}{c}{EO} & \multicolumn{2}{c}{GA} &
\multicolumn{2}{c}{heuristics} \\
\hline
{\em Hammond\/}    & 90   & (42s)  & 90   & (1s)    & 97  & (8s) \\
\multicolumn{7}{l}{\quad($N=4720$; $\alpha=5.8$)} \\
{\em Barth5\/}     & 139  & (64s)  & 139  & (44s)   & 146 & (28s) \\
\multicolumn{7}{l}{\quad($N=15606$; $\alpha=5.8$)} \\
{\em Brack2\/}     & 731  & (12s)  & 731  & (255s)  & \multicolumn{2}{c}{---} \\
\multicolumn{7}{l}{\quad($N=62632$; $\alpha=11.7$)} \\
{\em Ocean\/}      & 464  & (200s) & 464  & (1200s) & 499 & (38s) \\
\multicolumn{7}{l}{\quad($N=143437$; $\alpha=5.7$)} \\
\hline
\hline
Graph &&& \multicolumn{2}{c}{EO} & \multicolumn{2}{c}{SA} \\
\hline
{\em Nasa1824\/} &&& 739  & (6s)   & 739  & (3s)\\
\multicolumn{7}{l}{\quad($N=1824$; $\alpha=20.5$)} \\
{\em Nasa2146\/} &&& 870  & (10s)  & 870  & (2s)\\
\multicolumn{7}{l}{\quad($N=2146$; $\alpha=32.7$)} \\
{\em Nasa4704\/} &&& 1292 & (15s)  & 1292 & (13s)\\
\multicolumn{7}{l}{\quad($N=4704$; $\alpha=21.3$)} \\
{\em Stufe10\/}  &&& 51   & (180s) & 371  & (200s) \\
\multicolumn{7}{l}{\quad($N=24010$; $\alpha=3.8$)} \\
\end{tabular}
\end{center}
\label{tab1}
\end{table}

\section{Optimizing near Critical Points}

Further comparison of EO and SA, averaged over a large sample of a 
particular type of graph, shows EO to be especially 
useful near critical points (Boettcher, to appear).
It has been observed that many optimization problems exhibit critical
points delimiting ``easy'' phases of a generally hard
problem (Cheeseman 1991). Near such a critical point, finding solutions
becomes particularly difficult for local search methods that explore
some neighborhood in configuration space starting from an existing state.
Near-optimal solutions become widely
separated with diverging barrier heights between them. It is not
surprising that equilibrium search methods based on heat-bath techniques
like SA are not particularly successful here
(Binder 1987). In contrast, the driven dynamics of EO does not possess
any temperature control parameters that could increasingly limit the scale
of fluctuations.  A non-equilibrium approach like EO thus provides
a general-purpose optimization method that is complementary to SA,
which would be expected to freeze quickly into a poor local optimum
``where the {\em really\/} hard problems are'' (Cheeseman 1991).

\begin{figure}
\vskip 2.0truein
\includegraphics{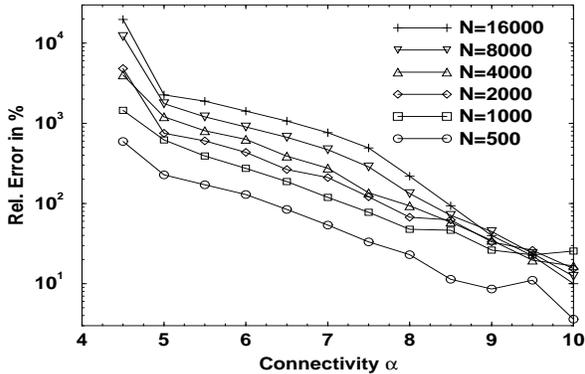}
\caption{Plot of SA's error relative to the best result found on 
geometric graphs, as a function of the mean connectivity $\alpha$. }
\label{error}
\end{figure}

\begin{figure}
\vskip 2.0truein
\includegraphics{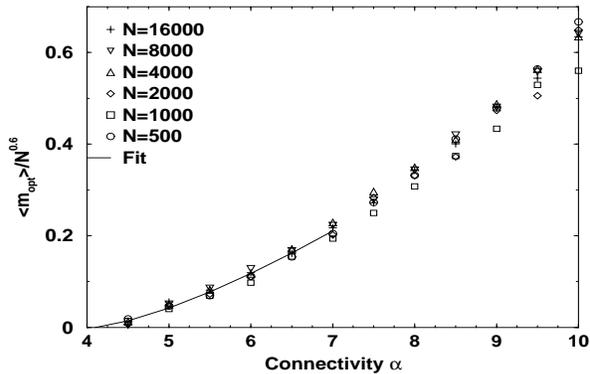}
\caption{Scaling plot of the data from EO according to Eq.~\protect\ref{scaleq}
for geometric graphs, as a function of the mean connectivity $\alpha$. The
scaling parameters and the fit are as discussed in the text.
}
\label{scaling}
\end{figure}

As an example, we explore this critical point for the equal
partitioning of geometric graphs, as a function of their
connectivity.\footnote{More
results of this study, including many different types of graphs, can be found
in (Boettcher, to appear).} It is
hopeless to obtain reliable benchmarks for the exact optimal partition
of large graphs.  Instead, by averaging over many
instances we can try to reproduce well-known results from the
percolation properties of this class of graphs. For instance, when the 
average connectivity $\alpha$ of a geometric graph is much below 
$\alpha_{\rm crit}\approx4.5$, the percolation threshold found for these
graphs (Balberg 1985), the graph most likely consists of many small clusters.
These can easily be sorted into equal sized partitions with vanishing cutsize,
at a cost of at most O($N^2$). When, on the other hand, the connectivity is
large, the graph is dense and almost homogeneous with many near-optimal
solutions in close 
proximity. But for connectivities near $\alpha_{\rm crit}$, a ``percolating''
cluster of size O($N$) appears with very widely separated minima (see 
Fig.~\ref{geograph}), making both the decision problem and the 
actual search very costly (Cheeseman 1991).

We have generated geometric graphs of connectivities between $\alpha=4$
and $\alpha=10$ (by varying the threshold distance $d$ below which points
are connected), at $N=500$, 1000, 2000, 4000, 8000, 
and 16000.  For each $\alpha$ we generated 16 different instances of 
graphs, identical for SA and EO. We performed 32 optimization runs for 
each method on each instance.  On each run, we used a different random seed to
establish an initial partition of the points.  SA was run using
the algorithm developed by Johnson {\em et al.\/} (Johnson 1989) for this case, 
but with a temperature length four times longer, to improve results.
EO was run for $200N$ update steps to produce a comparable runtime.
For each method, we have taken only the best result from all runs on a
given instance.  We average those best results,
for a particular connectivity $\alpha$, to obtain the mean cutsize for
that method as a function of $\alpha$ and $N$. To compare EO and SA,
we determine the relative error of SA with respect to the best result found
by either method (most often by EO!) for $\alpha\geq\alpha_{\rm crit}$.
Fig.~\ref{error} suggests that the error of SA diverges about linearly with
increasing $N$, near $\alpha_{\rm crit}$. 

For the data obtained with EO, we make an Ansatz
\begin{eqnarray}
\langle m_{\rm opt}\rangle\sim N^\nu\left(\alpha-\alpha_{\rm
0}\right)^\beta
\label{scaleq}
\end{eqnarray}
with $\nu=0.6$, in order to scale the data for all $N$ onto a single curve
(see Fig.~\ref{scaling}).  The remaining parameters are established according
to a data fit, yielding $\alpha_0=4.1$ and $\beta=1.4$. The fact that
$\alpha_0 < \alpha_{\rm crit}$ indicates that below the percolation threshold
EO's cutsizes are already non-vanishing, and so even EO does not always find
optimal partitions there.

\section{Extremal Optimization of the TSP}

In the graph partitioning problem, the implementation of EO is
particularly straightforward.  The concept of fitness, however, is
equally meaningful in any optimization problem whose cost function can
be decomposed into $N$ equivalent degrees of freedom.  Thus, EO may be
applied to many other NP-hard problems, even those where the choice of
quantities for the fitness function, as well as the choice of
elementary move, is less clear than in graph partitioning.  One case
where these choices are far from obvious is the traveling salesman
problem.  Even so, we have found there that EO presents a challenge
to more finely tuned methods (Boettcher, submitted).

In the traveling salesman problem (TSP), $N$ points (``cities'') are
given, and every pair of cities $i$ and $j$ is separated by a distance
$d_{ij}$. The problem is to connect the cities using the {\em
shortest\/} closed ``tour'', passing through each city exactly once.
For our purposes, take the $N\times N$ distance matrix $d_{ij}$ to be
symmetric. Its entries could be the Euclidean distances between cities
in a plane --- or alternatively, random numbers drawn from some
distribution, making the problem non-Euclidean.  (The former case
might correspond to a business traveler trying to minimize driving
time; the latter to a traveler trying to minimize expenses on a string
of airline flights, whose prices certainly do not obey triangle
inequalities!)

For the TSP, we implement EO in the following way.  Consider each city
$i$ as a degree of freedom, with a fitness based on the two links
emerging from it.  Ideally, a city would want to be connected to its
first and second nearest neighbor, but is often ``frustrated'' by the
competition of other cities, causing it to be connected instead to
(say) its $p$th and $q$th neighbors, $1\leq p\neq q\leq
N-1$.  Let us define the fitness of city $i$ to be
$\lambda_i=3/(p_i+q_i)$, so that $\lambda_i=1$ in the ideal
case.

Defining a move class (step (3) in EO's algorithm) is more difficult
for the TSP than for graph partitioning, since the constraint of a
closed tour requires an update procedure that changes several links at
once.  One possibility, used by SA among other local search methods,
is a ``two-change'' rearrangement of a pair of non-adjacent segments in
an existing tour.  There are $O(N^2)$ possible choices for a
two-change.  Most of these, however, lead to even worse results.
For EO, it would not be sufficient to select two independent cities of
poor fitness from the rank list, as the resulting two-change would
destroy more good links than it creates.  Instead, let us select one
city $i$ according to its fitness rank $n_i$, using the distribution
$P(n)\sim n^{-\tau}$ as before, and eliminate the longer of the two
links emerging from it.  Then, reconnect $i$ to a close neighbor,
using the {\em same\/} distribution function $P(n)$ as for the rank
list of fitnesses, but now applied instead to a rank list of $i$'s
neighbors ($n=1$ for first neighbor, $n=2$ for second neighbor, and so
on).   Finally, to form a valid closed tour, one of the old links to
the new (neighbor) city must be replaced; there is a unique way of doing so.
For the
optimal choice of $\tau$, this move class allows us the opportunity to
produce many good neighborhood connections, while maintaining enough
fluctuations to explore the configuration space.

\begin{table}[t]
\caption{Best tour-lengths found for the Euclidean (top) and the 
random-distance TSP (bottom).  Results for each value of $N$ are 
averaged over 10 instances, using on each instance an exact algorithm 
(except for $N=256$ Euclidean where none was available), the best-of-ten
EO runs, and the 
best-of-ten SA runs.  Euclidean tour-lengths are rescaled by $1/\sqrt{N}$.}
\begin{center}
\begin{tabular}{rrrr}
$N$ & Exact  & EO$_{10}$  & SA$_{10}$\\ 
\hline
Euclidean\qquad
 16& 0.71453&  0.71453& 0.71453\\
 32& 0.72185&  0.72237& 0.72185\\
 64& 0.72476&  0.72749& 0.72648\\
128& 0.72024&  0.72792& 0.72395\\
256& \qquad --- &  0.72707& 0.71854\\
\hline
Rand.\ Dist.\qquad
 16&  1.9368&  1.9368&  1.9368\\
 32&  2.1941&  2.1989&  2.1953\\
 64&  2.0771&  2.0915&  2.1656\\
128&  2.0097&  2.0728&  2.3451\\
256&  2.0625&  2.1912&  2.7803
\end{tabular}
\end{center}
\label{tab2}
\end{table}

We performed simulations at $N=16$, 32, 64, 128, and 256, in each case
generating ten random instances for both the Euclidean and
non-Euclidean TSP.  The Euclidean case consisted of $N$ points placed
at random in the unit square with periodic boundary conditions; the
non-Euclidean case consisted of a symmetric $N\times N$ distance
matrix with elements drawn randomly from a uniform distribution on the
unit interval.  On each instance we ran both EO and SA, selecting for
both methods the best of 10 runs from random initial conditions.  EO
used $\tau=4$ (Eucl.) and $\tau=4.4$ (non-Eucl.), with $16N^2$ update
steps.  SA used an annealing schedule with $\Delta T/T=0.9$ and
temperature length $32N^2$.  The results are given in
Table~\ref{tab2}, along with baseline results using an exact
algorithm.  While the EO results trail those of SA by up
to about 1\% in the Euclidean case, EO significantly outperforms SA
for the non-Euclidean (random distance) TSP.  Surprisingly, using
increased run times (longer temperature schedules) diminishes rather
than improves SA's performance in the latter case.  Finally, note that
one would not expect a general method such as EO to be competitive
here with specialized optimization algorithms designed particularly
with the TSP in mind.  But remarkably, EO's performance in both the
Euclidean and non-Euclidean cases --- within several percent of
optimality for $N\le 256$ --- places it not far behind the leading
specially-crafted TSP heuristics (Johnson 1997).

\section{Extremal Optimization and Learning}
Our results therefore indicate that a simple extremal optimization approach
based on self-organizing dynamics can outperform state-of-the-art
(and far more complicated or finely tuned) general-purpose algorithms on
hard optimization problems.
Based on its success on the generic and broadly applicable graph partitioning
problem, as well as on the TSP, we believe the concept will be applicable to
numerous other NP-hard problems.  It is worth stressing that the rank ordering
approach employed by EO is inherently non-equilibrium.  Such an approach
could not, for instance, be used to enhance SA, whose temperature schedule
requires equilibrium conditions.  This rank ordering serves as a sort of
``memory'', allowing EO to retain well-adapted pieces of a solution.  In this
respect it mirrors one of the crucial properties noted in the Bak-Sneppen
model (Boettcher 1996).  At the same time, EO maintains enough flexibility to
explore further reaches of the configuration space and to ``change its mind''.
Its success at this complex task provides motivation for the use of extremal
dynamics to model mechanisms such as learning, as has been suggested recently
to explain the high degree of adaptation observed in the brain (Chialvo 1999).

\subsubsection*{References}

C.~J.~Alpert and A.~B.~Kahng, 
{\em Integration: the VLSI Journal\/} {\bf 19}, 1 (1995).

P.~Bak, {\em How Nature Works\/} (Springer, New York, 1996).  

P.~Bak, C.~Tang, and K.~Wiesenfeld,
{\em Phys. Rev. Lett.\/} {\bf 59}, 381 (1987).

P.~Bak and K.~Sneppen, {\em Phys. Rev. Lett.\/} {\bf 71}, 4083
(1993).

I.~Balberg, Phys. Rev. B {\bf 31}, R4053 (1985).

K.~Binder, {\em Applications of the Monte Carlo
Method in Statistical Physics\/}, K.~Binder, Ed. (Springer, Berlin, 1987).

S.~Boettcher and M.~Paczuski, {\em Phys. Rev. E\/} {\bf 54}, 1082
(1996), and {\em Phys. Rev. Lett.\/} {\bf 79}, 889 (1997).

S.~Boettcher, {\em J.~Phys.~A:~Math.~Gen\/}, to appear; available
at http://xxx.lanl.gov/abs/cond-mat/9901353.

S.~Boettcher and A.~G.~Percus, submitted to \hfil\break
{\em Artificial~Intelligence\/};~available~at \hfil\break
http://xxx.lanl.gov/abs/cond-mat/9901351.

V.~{\v C}erny, {\em J.~Optimization Theory Appl.\/} {\bf 45}, 41 (1985).

P.~Cheeseman, B.~Kanefsky, and W.~M.~Taylor, in {\em
Proc. of IJCAI-91}, J.~Mylopoulos and R.~Rediter, Eds. (Morgan
Kaufmann, San Mateo, CA, 1991), pp.\ 331--337.

D.~R.~Chialvo and P.~Bak, {\em J. Neurosci.\/}, to appear.

C.~Darwin, {\em The Origin of Species by Means of Natural Selection\/}
(Murray, London, 1859).  

A.~E.~Dunlop and B.~W.~Kernighan, {\em IEEE Trans.
on Computer-Aided Design\/} {\bf CAD--4}, 92 (1985).

M.~R.~Garey and D.~S.~Johnson, {\em Computers
and Intractability: A Guide to the Theory of NP-Completeness\/} (Freeman,
New York, 1979).

S.~J.~Gould and N.~Eldridge, {\em Paleobiology\/} {\bf 3}, 115--151 (1977).

C. de Groot, D. Wuertz, K. H. Hoffmann,
{\em Lecture Notes in Computer Science\/} {\bf 496}, 445-454 (1991).

B.~A.~Hendrickson and R.~Leland, in {\em Proceedings of the 1995 
ACM/IEEE Supercomputing Conference (Supercomputing '95)\/}, San Diego,
CA, December 3--8, 1995 (ACM Press, New York, 1996).

J.~Holland, {\em Adaptation in Natural and Artificial Systems\/} (University 
of Michigan Press, Ann Arbor, 1975).

%C.~Hurwitz, {\em GNU tsp\_solve\/}, available at:
%\newline
%http://www.cs.sunysb.edu/\~{}algorith/implement/tsp/implement.shtml.

D.~S.~Johnson, C.~R.~Aragon, L.~A.~McGeoch, and 
C.~Schevon, {\em Operations Research\/} {\bf 37}, 865 (1989).  

D.~S.~Johnson and L.~A.~McGeoch, in {\em Local Search in Combinatorial
Optimization\/}, E.~H.~L.~Aarts and J.~K.~Lenstra, Eds. (Wiley, New York,
1997), chap.~8.

S.~Kirkpatrick, C.~D.~Gelatt, and M.~P.~Vecchi, {\em
Science\/} {\bf 220}, 671 (1983).

B. B. Mandelbrot, {\em The Fractal
Geometry of Nature\/} (Freeman, New York, 1983).

P.~Merz and B.~Freisleben, 
in {\em Lecture Notes in Computer Science: Parallel Problem Solving From
Nature --- PPSN V\/}, A.~E.~Eiben, T.~B\"ack, M.~Schoenauer, and H.-P.~Schwefel,
Eds. (Springer, Berlin, 1998), vol. 1498, pp.\ 765--774, and
{\em Technical Report No. TR--98--01\/}
(Department of Electrical Engineering and Computer Science, University of
Siegen, Siegen, Germany, 1998), available at
%\newline
http://www.informatik.uni-siegen.de/\~{}pmerz/publications.html.

N.~Metropolis, A.~W.~Rosenbluth, M.~N.~Rosenbluth,
A.~H.~Teller, and E.~Teller, {\em J.~Chem. Phys.\/} {\bf 21}, 1087 (1953).

M.~Paczuski, S.~Maslov, and
P.~Bak, {\em Phys. Rev. E\/} {\bf 53}, 414 (1996).  

I. Rodriguez-Iturbe and A. Rinaldo, {\em
Fractal River Basins: Chance and Self-Organization\/} (Cambridge, New
York, 1997).

K.~Sneppen, P.~Bak, H.~Flyvbjerg, and M.~H.~Jensen,
{\em Proc. Natl. Acad. Sci.\/} {\bf 92}, 5209 (1995).

\end{document}